\documentclass{tac}

\usepackage{amsmath}
\usepackage{amssymb}
\usepackage{xy}

\input diagxy

\usepackage[colorlinks=true]{hyperref}
\hypersetup{allcolors=[rgb]{0.1,0.1,0.4}}

\author{David Neal Broodryk}

\address{Department of Mathematics and Applied Mathematics, University of Cape Town\\
Rondebosch 7701\\[5pt]
}

\title{Characterization of coextensive varieties of universal algebras}

\copyrightyear{2020}


\eaddress{BRDDAV020@myuct.ac.za}

\newtheorem{thm}{Theorem}

\newtheoremrm{rem}{Remark}

\let\pf\proof
\let\epf\endproof

\newcommand{\Fxy}{F(\{x,y\})}
\newcommand{\C}{\mathcal{C}}
\newcommand{\Fo}{F(\emptyset)}
\renewcommand{\thesubsection}{\arabic{subsection}}

\begin{document}

\maketitle
\begin{abstract}
A coextensive category can be defined as a category $\C$ with finite products such that for each pair $X,Y$ of objects in $\C$, the canonical functor $\times\colon X/\C \times Y/\C \to (X \times Y)/\C$ is an equivalence. We give a syntactical characterization of coextensive varieties of universal algebras.
\end{abstract}

An extensive category can be defined as a category $\C$ with finite coproducts such that for each pair $X,Y$ of objects in $\C$, the canonical functor $+\colon \C/X \times \C/Y \to \C/(X + Y)$ is an equivalence. A category that satisfies the dual condition is called coextensive. According to \cite{CLW}, the term ``extensive category'' was first used by W.\ F.\ Lawvere and S.\ Schanuel, although ``categories with disjoint and universal coproducts'' were considered by A.\ Grothendieck a long time ago, and there are related papers of various authors. Examples of extensive categories include the category $\mathbf{Cat}$ of all small categories and the category $\mathbf{Top}$ of topological spaces. An example of a coextensive category is the category $\mathbf{CRing}$ of commutative rings.

\definition\label{definition}  

A category $\C$ with finite products is called coextensive if for each pair $X,Y$ of objects in $\C$, the canonical functor $\times\colon X/\C \times Y/\C \to (X \times Y)/\C$ is an equivalence. Equivalently let $L$ be the left adjoint of $\times$ and let $\mu$ and $\varepsilon$ be the unit and counit of this adjunction, then $L\colon  (X \times Y)/\C \to X/\C \times Y/\C$ sends $v \colon  (X \times Y) \to Z$ to the pair of canonical maps $(i_1,i_2)\colon  (X,Y) \to (X +_{X \times Y} Z,Y +_{X \times Y} Z)$ and $\C$ is coextensive if and only if both $\varepsilon$ and $\mu$ are natural isomorphisms.
\endthm

As follows from Proposition 2.2 of \cite{CLW}, $\C$ is coextensive if and only if for for every commutative diagram
$$\bfig
 \square[X \times Y`Y`A`C;\pi_2`h`g`p_2]
 \morphism(0,0)|b|<-500,0>[A`B;p_1]
 \morphism(0,500)|a|<-500,0>[X \times Y`X;\pi_1]
 \morphism(-500,500)|a|<0,-500>[X`B;f]
\efig$$

The two squares are pushouts if and only if the bottom row is a product diagram. Let $\C$ be a variety of universal algebras and let $F(X)$ denote the free algebra in $\C$ on the set $X$. In particular we shall use $\Fxy$ which is the algebra consisting of terms of at most two variables. We shall also use $\Fo$ which is the algebra consisting of all constant terms and is the initial object in $\C$.

\proposition\label{proposition}  Any coextensive variety contains a $(k+2)$-ary term $t$ and constants $e_1, \dots,e_k,e_1', \dots,e_k' \in \Fo$ such that the identities hold:

 \begin{align*}
    t(x, y, e_1, \dots, e_k) &= x \\
    t(x, y, e_1', \dots, e_k') &= y
\end{align*}

\endthm

\pf Let $A\subseteq \Fxy \times \Fxy$ be the subalgebra generated by $\Fo \times \Fo \cup \{(x,x),(y,y)\}$. Consider the diagram where $i:A \to \Fxy \times \Fxy$ is the inclusion:

$$\bfig
 \morphism(0,1000)<0,-1000>[\Fo \times \Fo`A;\subset]
 \morphism(0,1000)<1000,0>[\Fo \times \Fo`\Fo;\pi_2]
 \morphism(0,0)<1000,0>[A`\Fxy; p_2\circ i]
 \morphism(0,0)|a|<-1000,0>[A`\Fxy; p_1\circ i]
 \morphism(0,1000)|a|<-1000,0>[\Fo \times \Fo`\Fo;\pi_1]
 \morphism(-1000,1000)|a|<0,-1000>[\Fo`\Fxy;\subset]
 \morphism(1000,1000)|a|<0,-1000>[\Fo`\Fxy;\subset]
 \morphism(0,0)<0,-1000>[A`\Fxy \times \Fxy;i]
 \morphism(-1000,0)<0,-1000>[\Fxy`\Fxy;=]
 \morphism(1000,0)<0,-1000>[\Fxy`\Fxy;=]
 \morphism(0,-1000)<-1000,0>[\Fxy \times \Fxy`\Fxy;p_1]
 \morphism(0,-1000)<1000,0>[\Fxy \times \Fxy`\Fxy;p_2]
\efig$$

Note that since $A$ contains the diagonal $\Delta(\Fxy)$ the maps $p_1\circ i$ and $p_1\circ i$ are surjective. On the other hand the elements of $A$ are all of the form $t((x,x),(y,y),(e_1,e_1'),\dots,(e_k,e_k'))$ for some term $t$ and constants $e_i,e_i'$. Therefore the quotients of $A$ induced by $p_1\circ i$ and $p_2\circ i$ are exactly the quotients induced by $\pi_1$ and $\pi_2$ respectively. We conclude that the top squares are pushout diagrams and therefore that the middle row is a product diagram. In particular $(x,y) \in A$, which implies:

\[(x,y) = t((x,x),(y,y),(e_1,e_1'),\dots,(e_k,e_k')) = (t(x,y,e_1,\dots,e_k),t(x,y,e_1',\dots,e_k'))\]

for some term $t$ and constants $e_i,e_i'$ as desired.

\epf

We call such an operation t a diagonalizing term and in any algebra $A \supseteq \Fo \times \Fo$ we define the map $\delta(x,y) = t(x,y,(e_1,e_1'),\dots,(e_k,e_k'))$

\proposition\label{proposition_3}
	Let $\C$ be a variety of universal algebras containing a diagonalizing term $t$ equipped with constants $e_1,\dots, e_k, e_1', \dots, e_k'$. Then the following conditions hold:

	\begin{enumerate}
		\item $\C$ is left coextensive
		\item for any $X,Y \in \C$ and any surjective $q: X\times Y \to Z$ there exist surjective $f: X \to X'$ and $g: Y \to Y'$ such that $q \simeq f\times g$.
	\end{enumerate}

	\pf 
	$\C$ is left coextensive by Theorem 3 of \cite{DB} since we can fix any constant as $0$ and have the equalities:
	\begin{align*}
	    t(x, 0, e_1, \dots, e_k) &= x \\
	    t(x, 0, e_1', \dots, e_k') &= 0
	\end{align*}

	\vspace{1em}

	Let $q: X\times Y \to Z$ be surjective. We define the relations:

	\begin{itemize}
		\item $E_X = \{(a,c)\in X\times X\hspace{0.4em} | \hspace{0.4em} \exists_{b,d \in Y} : q((a,b))=q((c,d))\}$
		\item $E_Y = \{(b,d)\in Y\hspace{0.15em}\times Y \hspace{0.55em} | \hspace{0.4em} \exists_{a,c \in X} : q((a,b))=q((c,d))\}$
	\end{itemize}

	In any variety these relations are reflexive, symmetric and homomorphic. To show that they are also transitive let $(a,c),(c,f) \in E_X$. Then there exist $b,d,e,g$ such that $q((a,b))=q((c,d))$ and $q((c,e))=q((f,g))$.

	\begin{align*}
	    q((a,e)) &= q(t( (a,b), (c,e), (e_1,e_1'), \dots, (e_k,e_k') ))\\ 
	    &= q(t( (c,d), (c,e), (e_1,e_1'), \dots, (e_k,e_k') )) \\
	    &= q((c,e)) \\
	    &= q((f,g))
	\end{align*}

	Therefore $E_X$ and $E_Y$ are congruences and we can define the quotient projections: $q_X: X \to X/E_X$, $q_Y: Y \to Y/E_Y$, and $q_x \times q_y: X \times Y \to (X\times Y)/(E_x \times E_y)$. Clearly if $q((a,b))=q((c,d))$ then $(q_x \times q_y)((a,b))=(q_x \times q_y)((c,d))$. On the other hand, if $(q_x \times q_y)((a,b))=(q_x \times q_y)((c,d))$ then there exist $a',c'\in Y$ and $b',d' \in X$ such that $q((a,a')) = q((c,c'))$ and $q((b',b)) = q((d',d))$ Therefore we have:

	\begin{align*}
	    q(a,b) &= q(t( (a,a'), (b',b), (e_1,e_1'), \dots, (e_k,e_k') ))\\ 
	    &= q(t( (c,c'), (d',d), (e_1,e_1'), \dots, (e_k,e_k') )) \\
	    &= q(c,d)
	\end{align*}

	So $q \simeq q_x \times q_y$ as desired.

	\epf

\endthm

\proposition\label{proposition_3}  

Let $\C$ be a variety, the following statements are equivalent:

\begin{enumerate}
	\item $\C$ is coextensive
	\item $\C$ contains a diagonalizing term and for any $A\in \C$ and any equalities in $FU(A)$: 
	\begin{align*}
	    u(a_1, \dots, a_m, \omega_1, \dots, \omega_n) &= v(a_1, \dots, a_m, \omega_1, \dots, \omega_n) \\
	    u(a_1, \dots, a_m, \omega_1', \dots, \omega_n') &= v(a_1, \dots, a_m, \omega_1', \dots, \omega_n')
	\end{align*}

	we have the equality in $FU(A) + \Fo^2$: 

	\[ u(a_1, \dots, a_m, (\omega_1,\omega_1'), \dots, (\omega_n,\omega_1')) = v(a_1, \dots, a_m, (\omega_1,\omega_1'), \dots, (\omega_n,\omega_1')) \]

	\item $\C$ contains a diagonalizing term and for any set $X$ satisfies the equations:
	\begin{enumerate}
	\item $\delta(x,x) = x$ for any $x\in F(X) + \Fo^2$
	\item $u(\delta(x_1,y_1),\dots,\delta(x_n,y_n)) = \delta(u(x_1,\dots,x_n),u(y_1,\dots,y_n))$ for any operation $u$ and any $x_i,y_i \in F(X) + \Fo^2$
	\end{enumerate}

\end{enumerate}

\endthm

\pf $(1 \iff 2)$ If $\C$ is coextensive then by Proposition 2, $\C$ contains a diagonalizing term $t$. On the other hand if $\C$ contains such a term then by Proposition 3.1, $\C$ is left coextensive. It remains to show that in any such $\C$, when diagram 1 is a pushout diagram it is also a product diagram. It is sufficient to consider only $X = Y = \Fo$. In particular consider for any algebra $A$ equipped with a morphism $h: \Fo^2 \to A$ the following pushout diagram:

$$\bfig
 \morphism(0,1000)<0,-1000>[\Fo \times \Fo`FU(A) + \Fo^2;\subset]
 \morphism(0,1000)<1000,0>[\Fo \times \Fo`\Fo;\pi_2]
 \morphism(0,0)<-1000,0>[FU(A) + \Fo^2`FU(A); p_1]
 \morphism(0,0)<1000,0>[FU(A) + \Fo^2`FU(A); p_2]
 \morphism(0,1000)|a|<-1000,0>[\Fo \times \Fo`\Fo;\pi_1]
 \morphism(-1000,1000)|a|<0,-1000>[\Fo`FU(A);\subset]
 \morphism(1000,1000)|a|<0,-1000>[\Fo`FU(A);\subset]
 \morphism(0,0)<0,-1000>[FU(A) + \Fo^2`A;{[\varepsilon_A, h]}]
 \morphism(-1000,0)<0,-1000>[FU(A)`B;f]
 \morphism(1000,0)<0,-1000>[FU(A)`C;g]
 \morphism(0,-1000)<-1000,0>[A`B;p'_1]
 \morphism(0,-1000)<1000,0>[A`C;p'_2]
\efig$$

Where $U$ is the underlying set functor and $+$ denotes the coproduct. All four small squares are defined to be pushouts. In particular note that $[\varepsilon_A,h]$ is surjective, so if the first row is a product diagram, then by Proposition 3.2 so is the second row and the whole square. Therefore it is sufficient to consider only the top row.

\vspace{1em}

Consider the map $(p_1, p_2) : FU(A) + \Fo^2 \to FU(A)^2$. $\C$ will be coextensive if and only if this map is a bijection. Note that for any $a,b \in FU(A)$ we have:

\[ (p_1, p_2)(\delta(a,b)) = (t(a,b,e_1,\dots,e_k), t(a,b,e_1',\dots,e_k')) = (a,b) \]

and so $(p_1, p_2)$ is surjective. On the other hand it is injective if and only if for all pairs of identities:

 \begin{align*}
    u(a_1, \dots, a_m, \omega_1, \dots, \omega_n) &= v(a_1, \dots, a_m, \omega_1, \dots, \omega_n) \\
    u(a_1, \dots, a_m, \omega_1', \dots, \omega_n') &= v(a_1, \dots, a_m, \omega_1', \dots, \omega_n')
\end{align*}

we have: 

\[ u(a_1, \dots, a_m, (\omega_1,\omega_1'), \dots, (\omega_n,\omega_1')) = v(a_1, \dots, a_m, (\omega_1,\omega_1'), \dots, (\omega_n,\omega_1')) \]

for some operations $u,v$, variables $a_1,\dots,a_m$, and constants $\omega_1,\dots,\omega_n,\omega_1',\dots,\omega_n'$.

\vspace{1em}

$(2 \iff 3)$ Let $X$ be any set and note that if $2$ holds we immediately have:

\begin{enumerate}
\item $\delta(x,x) = x$ for any $x\in F(X) + \Fo^2$
\item $u(\delta(x_1,y_1),\dots,\delta(x_n,y_n)) = \delta(u(x_1,\dots,x_n),u(y_1,\dots,y_n))$ for any operation $u$ and any $x_i,y_i \in F(X) + \Fo^2$
\end{enumerate}

On the other hand if these two equations hold then for any such pair of identities we have:

 \begin{align*}
    u(a_1, \dots, a_m, (\omega_1,\omega_1'), \dots, (\omega_n,\omega_1')) &= u(\delta(a_1,a_1), \dots, \delta(a_m,a_m), \delta(\omega_1,\omega_1'), \dots, \delta(\omega_n,\omega_n')) \\
    &= \delta(u(a_1, \dots, a_m, \omega_1, \dots, \omega_n), u(a_1, \dots, a_m, \omega_1', \dots, \omega_n')) \\
    &= \delta(v(a_1, \dots, a_m, \omega_1, \dots, \omega_n), v(a_1, \dots, a_m, \omega_1', \dots, \omega_n')) \\
    &= v(\delta(a_1,a_1), \dots, \delta(a_m,a_m), \delta(\omega_1,\omega_1'), \dots, \delta(\omega_n,\omega_n')) \\
    &=v(a_1, \dots, a_m, (\omega_1,\omega_1'), \dots, (\omega_n,\omega_1'))
\end{align*}

\epf

\thm\label{theorem} 
	A variety $\C$ is coextensive if and only if the following conditions hold: 

	\begin{enumerate}

		\item $\C$ contains a $(k+2)$-ary term $t$ and constants $e_1, \dots,e_k,e_1', \dots,e_k' \in \Fo$ such that the identities hold:

		\begin{align*}
		    t(x, y, e_1, \dots, e_k) &= x \\
		    t(x, y, e_1', \dots, e_k') &= y
		\end{align*}
		
		\item There exist terms $\alpha_0,\dots,\alpha_n,\beta_0,\dots,\beta_n \in F(\{x\} + U(\Fo^2))$, operations $u_0,\dots,u_m \in F(X)$ and for all operations $s \in F(X)$ there exist $u^{(s)}_0,\dots,u^{(s)}_m$ such that

		\begin{align*}
		    u_0(\beta_0,\dots,\beta_n) &= \delta(x,x) \\
		    u_m(\alpha_0,\dots,\alpha_n) &= x \\
		    u^{(s)}_0(\beta_0,\dots,\beta_n) &= \delta(s(x_1,\dots,x_l),s(x_1',\dots,x_l')) \\
		    u^{(s)}_m(\alpha_0,\dots,\alpha_n) &= s(\delta(x_1,x_1'),\dots,\delta(x_l,x_l'))
		\end{align*}

		and for all $i<m$:

		\begin{align*}
		    u_i(\alpha_0,\dots,\alpha_n) &= u_i(\beta_0,\dots,\beta_n) \\
		    u^{(s)}_i(\alpha_0,\dots,\alpha_n) &= u^{(s)}_i(\beta_0,\dots,\beta_n) \\
		\end{align*}

		\item For each $j<n$ one of the following is true:

		\begin{enumerate}
			\item $\alpha_j = \beta_j$
			\item $\alpha_j = (\omega,\omega')$ and $\beta_j = \delta(\omega,\omega')$ for some constants $\omega,\omega'$
			\item $\beta_j = (\omega,\omega')$ and $\alpha_j = \delta(\omega,\omega')$ for some constants $\omega,\omega'$
			\item $\alpha_j = v(\delta(\omega_0,\omega_0'),\dots,\delta(\omega_l,\omega_l'))$ and $\beta_j = \delta(v(\omega_0,\dots,\omega_l),v(\omega_0',\dots,\omega_l'))$ for some constants $\omega_0,\dots,\omega_l,\omega_0,\dots,\omega_l'$
			\item $\beta_j = v(\delta(\omega_0,\omega_0'),\dots,\delta(\omega_l,\omega_l'))$ and $\alpha_j = \delta(v(\omega_0,\dots,\omega_l),v(\omega_0',\dots,\omega_l'))$ for some constants $\omega_0,\dots,\omega_l,\omega_0,\dots,\omega_l'$
		\end{enumerate}

	\end{enumerate}

where $delta(a,b) = t(a,b,(e_0,e_0'),\dots,(e_k,e_k'))$ and $X$ is any set.

\endthm

\pf Condition 1 simply says that $\C$ has a diagonalizing term. Let $1_F(X) + \varepsilon_{\Fo}: F(X) + FU(\Fo^2) \to F(X) + \Fo^2$ be the map induced by the identity at $F(X)$ and the counit at $\Fo$. Then $F(X) + \Fo^2 \simeq F(X \cup U(\Fo^2))/E$ where $E$ is the congruence generated by this map. Note that this congruence is generated by the relation R = $\{\delta(\omega,\omega') = (\omega,\omega')\} \cup \{u(\delta(\omega_1,\omega_1'),\dots,\delta(\omega_n,\omega_n')) = \delta(u(\omega_1,\dots,\omega_n),u(\omega_1',\dots,\omega_n'))\}$

Note that $\delta(x,x) = x$ for any $x\in F(X) + \Fo^2$ if and only if $(\delta(x,x),x) \in E$ if and only if there exists a natural number $n$ such that $(\delta(x,x),x) \in Q^n$, where $Q$ is the reflexive symmetric homomorphic relation on $F(X) + FU(\Fo^2)$ generated by $R$. This is true if and only if there exist $a_0, \dots, a_k \in F(X) + FU(\Fo^2)$ such that $a_0 = \delta(x,x)$, $a_k = x$ and $(a_i,a_{i+1}) \in Q$ for $i<k$. But $(a_i,a_{i+1}) \in Q$ if and only if for some term $u_i$ we have:

\begin{align*}
    a_i = u_i(\beta_0,\dots,\beta_n) \\ 
    a_{i+1} = u_i(\alpha_0,\dots,\alpha_n) \\ 
\end{align*}

Where the $\alpha_j$ and $\beta_j$ are in the symmetric reflexive relation generated by $R$. This is true exactly when condition 3 above is satisfied. Similarly for any operation $s$ and variables $x_1,\dots,x_n,x_1',\dots,x_n' \in X$ there exist $a_0', \dots, a_k' \in F(X) + FU(\Fo^2)$ such that $a_0' = s(\delta(x_1,x_1'),\dots,\delta(x_n,x_n'))$, $a_k' = \delta(s(x_1,\dots,x_n),s(x_1',\dots,x_n'))$ and for $i<k$ there exists some term $u_i'$ such that:

\begin{align*}
    a_i = u_i^{(s)}(\beta_0,\dots,\beta_n) \\ 
    a_{i+1} = u_i^{(s)}(\alpha_0,\dots,\alpha_n) \\ 
\end{align*}

Where the $\alpha_j$ and $\beta_j$ are in the symmetric reflexive relation generated by $R$. To simplify the notation we can assume that the $\alpha_j$ and $\beta_j$ are the same for every $u_i$ and $u_i^{(s)}$. Now by proposition 4, we have that $\C$ is coextensive if and only if the above conditions hold. 

\epf

Recall that $\delta(x,y)$ is not an operation in the usual sense but is instead defined as $\delta(x,y) = t(x,y,(e_1,e_1'),\dots,(e_k,e_k'))$ Therefore equations involving $\delta$ are actually equations involving $t$ and the elements $(e_1,e_1'),\dots,(e_k,e_k')$. Note that $F(X) + FU(\Fo^2) = F(X + U(\Fo^2))$ is a free algebra in which the $(\omega,\omega')$ are elements of the generating set. Therefore the equations of condition 2 in the theorem must still hold after substituting $(\omega,\omega') = x_{\omega,\omega'}$ for some variables $x_{\omega,\omega'}\in X$. As an example it is helpful to consider the variety of commutative rings.

\example\label{example}  

Let $\C$ be the variety of commutative rings. let $t(a,b,c,d) = a \cdot c + b \cdot d$, $e_0 = 1, e_1 = 0, e_0' = 0, e_1' = 1$. Then we have $t(a,b,e_0,e_1) = a\cdot1 + b\cdot 0 = a$ and $t(a,b,e_0',e_1') = a\cdot0 + b\cdot1 = b$, so $t$ is a diagonalizing term.

\vspace{1em}

To show $\delta(x,x) = x$ consider:

\begin{align*}
	u_0(a_0,a_1,a_2,a_3) &= a_0\cdot a_1 + a_0\cdot a_2 \\
	u_1(a_0,a_1,a_2,a_3) &= a_0\cdot a_3 \\
	\alpha &= \big(x,(1,0),(0,1),1\big) \\
	\beta &= \big(x,(1,0),(0,1),1\cdot(1,0) +1\cdot(0,1)\big)
\end{align*}

 Note in particular that $\alpha_3 = 1$ and $\beta_3 = 1\cdot(1,0) +1\cdot(0,1)$ satisfies condition $3(d)$ of the theorem. This gives the equations:

 \[\delta(x,x) = x\cdot(1,0) + x\cdot(0,1) = x\cdot(1\cdot(1,0) + 1\cdot(0,1)) \simeq_E x\cdot1 = x\]

Therefore $\delta(x,x) = x$ as desired. It remains to show that $s(\delta(x_1,y_1),\dots,\delta(x_n,y_n)) = \delta(s(x_1,\dots,x_n),u(y_1,\dots,y_n))$ for any operation $s$ and any $x_i,y_i \in FU(A) + \Fo^2$. It suffices to show this for $s(x,y) = x+y$ and $s(x,y)=x\cdot y$.

To show $\delta(x+x',y+y') = \delta(x,y)+\delta(x',y')$ consider:

\begin{align*}
	u_0(a_0,a_1,a_2,a_3,a_4,a_5) &= a_0\cdot a_4 + a_1\cdot a_4 + a_2\cdot a_5 + a_3\cdot a_5 \\
	\alpha &= \big( x,x',y,y',(1,0),(0,1) \big) \\
	\beta &=  \big( x,x',y,y',(1,0),(0,1) \big)
\end{align*}

 Note that in this case no equivalence under $E$ is necessary so we have $\alpha=\beta$. This gives the equations:

 \[\delta(x+x',y+y') = (x+x')\cdot(1,0) + (y+y')\cdot(0,1) = x\cdot(1,0) + y\cdot(0,1) + x'\cdot(1,0) + y'\cdot(0,1) = \delta(x,y) + \delta(x',y')\]

 To show $\delta(x\cdot x',y\cdot y') = \delta(x,y)\cdot \delta(x',y')$ consider:

\begin{align*}
	u_0(a_0,\dots,a_11) &= a_0a_1a_4 + a_0a_3a_6 + a_1a_2a_6 + a_2a_3a_5 \\
	u_1(a_0,\dots,a_11) &= a_0a_1a_7 + a_0a_3a_8 + a_1a_2a_8 + a_2a_3a_9 \\
	u_2(a_0,\dots,a_11) &= (a_0a_10 + a_2a_11)\cdot(a_1a_10 + a_3a_11) \\
	\alpha &= \big( x,x',y,y',(1,0),(0,1),0, \delta(1\cdot1,0\cdot0), \delta(1\cdot0,0\cdot1), \delta(0\cdot0,1\cdot1),\delta(1,0),\delta(0,1) \big) \\
	\beta &=  \big( x,x',y,y',\delta(1,0),\delta(0,1),\delta(0,0),\delta(1,0)^2,\delta(1,0)\delta(0,1),\delta(0,1)^2,(1,0),(0,1) \big)
\end{align*}

 This gives the equations:

\begin{align*}
	\delta(x\cdot x',y\cdot y') &= x x'\cdot (1,0) + y y'\cdot (0,1) \\
	&= x x'\cdot (1,0) + x y'\cdot 0 + y x'\cdot 0 + y y'\cdot (0,1)\\
	&\simeq_E x x'\cdot \delta(1,0) + x y'\cdot \delta(0,0) + y x'\cdot \delta(0,0) + y y'\cdot \delta(0,1) \\
	&= x x'\cdot \delta(1\cdot1,0\cdot0) + x y'\cdot \delta(1\cdot0,0\cdot1) + y x'\cdot \delta(1\cdot0,0\cdot1) + y y'\cdot \delta(0\cdot0,1\cdot1) \\
	&\simeq_E x x'\cdot \delta(1,0)^2 + x y'\cdot \delta(1,0)\delta(0,1) + y x'\cdot \delta(1,0)\delta(0,1) + y y'\cdot \delta(0,1)^2 \\
	&= (x\cdot \delta(1,0) + y\cdot \delta(0,1)) \cdot (x'\cdot \delta(1,0) + y'\cdot \delta(0,1)) \\
	&\simeq_E (x\cdot (1,0) + y\cdot (0,1)) \cdot (x'\cdot (1,0) + y'\cdot (0,1)) \\
	&= \delta(x,y)\cdot \delta(x',y')
\end{align*}

And so the variety of commutative rings is coextensive.

\endthm

\refs

\bibitem [1]{CLW} [1] A. Carboni, S. Lack, and R. F. C. Walters, Introduction to extensive and distributive categories, Journal of Pure and Applied Algebra, 84(2), 1993, 145-158.

\bibitem [2]{DB} [2] D. Broodryk, Characterization of left coextensive varieties of universal algebras, Theory and Applications of Categories, Vol. 34, No. 32, 2019, 1036-1038.

\endrefs

\end{document}